\documentclass{article}
\usepackage{arxiv}

\usepackage[utf8]{inputenc} 
\usepackage[T1]{fontenc}    
\usepackage{hyperref}       
\usepackage{url}            
\usepackage{booktabs}       
\usepackage{amsfonts}       
\usepackage{nicefrac}       
\usepackage{microtype}      
\usepackage{lipsum}
\usepackage{graphicx,xcolor}
\usepackage[fleqn]{amsmath}
\usepackage{newtxtext}
\usepackage[varg]{newtxmath}
\usepackage{latexsym}
\usepackage{newtxtext}
\usepackage[varg]{newtxmath}
\usepackage{color}
\usepackage{algorithm,algorithmic}
\usepackage{url}

\newcommand{\transpose}[1]{{#1}^{\top}}
\newtheorem{theorem}{Theorem}

\newcommand{\QED}{\hfill \ensuremath{\Box}}

\title{Fairness in Robust Unit Commitment Problem Considering Suppression of Renewable Energy}

\author{
 Ichiro Toyoshima \\
  Toshiba Energy Systems \& Solutions Corporation\\
  Fuchu, Tokyo, Japan \\
  \texttt{ichiro.toyoshima.r24@mail.toshiba} \\
   \And
 Pierre-Louis POIRION \\
  Center for Advanced Intelligence Project \\
  RIKEN \\
  Chuo, Tokyo, Japan \\
  \And
 Tomohide Yamazaki \\
  Toshiba Energy Systems \& Solutions Corporation\\
  Fuchu, Tokyo, Japan \\
\And
 Kota Yaguchi \\
  Toshiba Energy Systems \& Solutions Corporation\\
  Fuchu, Tokyo, Japan \\
\And
 Masayuki Kubota \\
  Toshiba Energy Systems \& Solutions Corporation\\
  Fuchu, Tokyo, Japan  \\
\And
 Ryota Mizutani \\
  Toshiba Energy Systems \& Solutions Corporation\\
  Kawasaki, Kanagawa, Japan \\
\And
 Akiko Takeda \\
 Graduate School of Information Science and Technology \\
 The University of Tokyo \\
Bunkyo, Tokyo, Japan \\
}

\begin{document}
\maketitle
\begin{abstract}
Power company operators make power generation plans one day in advance, in what is known as the Unit Commitment (UC) problem. UC is exposed to uncertainties, such as unknown electricity load and disturbances caused by renewable energy sources, especially PVs. In previous research, we proposed the Renewable Energy Robust Optimization Problem (RE-RP), which solves these uncertainties by considering suppression. In this paper, we propose a new model called RE-RP with fairness (RE-RPfair), which aims to achieve fair allocation among PVs allocation. This model is an expansion of the original RE-RP, and we prove its effectiveness through simulation. To measure the degree of fairness, we use the Gini Index, which is well-known in social science.
\end{abstract}


\section{Introduction}
\begin{figure}[bh]
\begin{center}
 \includegraphics[width=8.0cm]{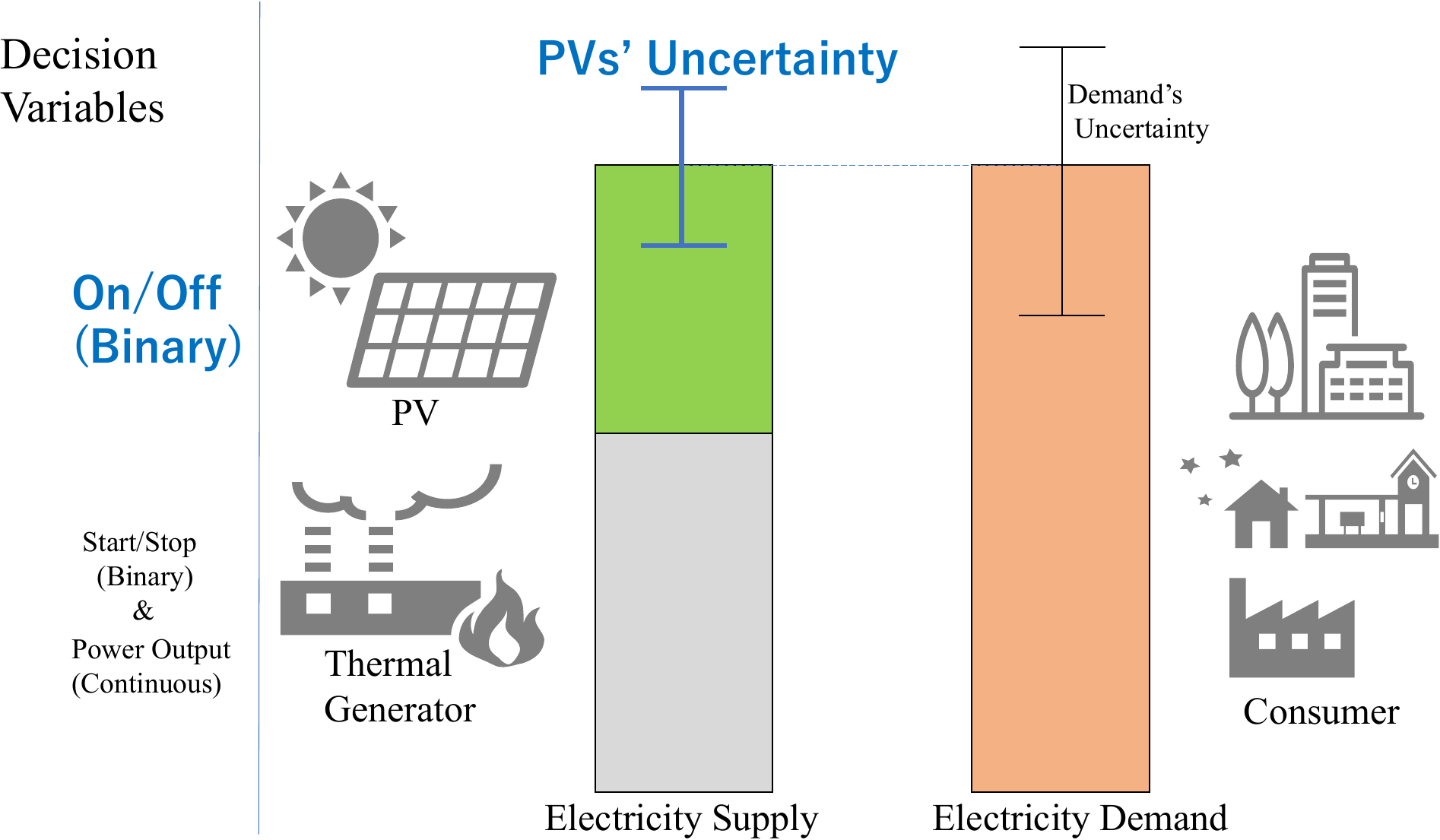}
\end{center}
\caption{Setting of RE-RP \ (On/Off binary variables and PVs' uncertainty are key features.)}\label{fig:Settings}
\end{figure}
Power company operators typically make power generation plans a day in advance with the intention of minimizing fuel costs. The schedule is obtained by solving the so-called \textit{Unit Commitment (UC)} problem, which is traditional topic in the energy management research field. Many prominent approaches, such as \cite{UCsurvey_old,UCsurvey_new}, have been proposed to solve conventional UC problems.

Recently, the UC problems are being exposed to uncertainties, which stem from two main sources. The first source is the unknown load for the next day, while the second source comes from renewable energy (RE) sources, particularly photovoltaic (PV) systems. However, even the most reliable forecasting techniques are not entirely error-free. Thus, models designed for UC problems must be equipped with mechanisms to handle uncertainties. Fortunately, many studies have been conducted on how to deal with uncertainty in energy systems, which can be found in the following references~\cite{takano2016,aihara2015,udagawa2016,murakami2011,shiina2003,zheng}.

However, we are currently facing a new challenge in addressing the UC problem, as another factor needs to be taken into account. PVs have spread more rapidly than anticipated, and as a result, system operators are often required to decide which PV output to suppress to maintain system security, or its ability to remain stable in response to disturbances, while also addressing UC issues. The suppression of PV output significantly impacts the uncertainties involved. Unfortunately, existing methods are not equipped to handle such uncertainties as their description varies depending on the decisions made. Thus, there is a need for a new method that can effectively address such uncertainties.

Prior to this work, we focused on a technique called \textit{robust optimization}~\cite{roboptbook}  as a way of handling the uncertainty in UC problems.
We note that some applications of robust optimization to energy systems existed~\cite{Bertsimas,Alvaro,moretti,zhaoguan,baziar,udagawa2017control}.
Unlike these papers, we needed to, in our setting, take into account both the uncertainties and  each PV unit suppression.
More specifically, in order to suppress a PV unit, an ``on/off'' signal is typically used.

Our previous work~\cite{CIGREyamazaki}
proposed the Renewable Energy Robust optimization Problem (RE-RP), a optimization problem that expands the model in~\cite{Bertsimas} which does not take into account PV power.
Fig.~\ref{fig:Settings} shows the concept of RE-RP.
They demonstrated the effectiveness of RE-RP through experiments that simulated a power company with several decades of thermal units and PVs located in~\cite{CIGREyamazaki}. 

\begin{figure}[tbh]
\begin{center}
 \includegraphics[width=9.0cm]{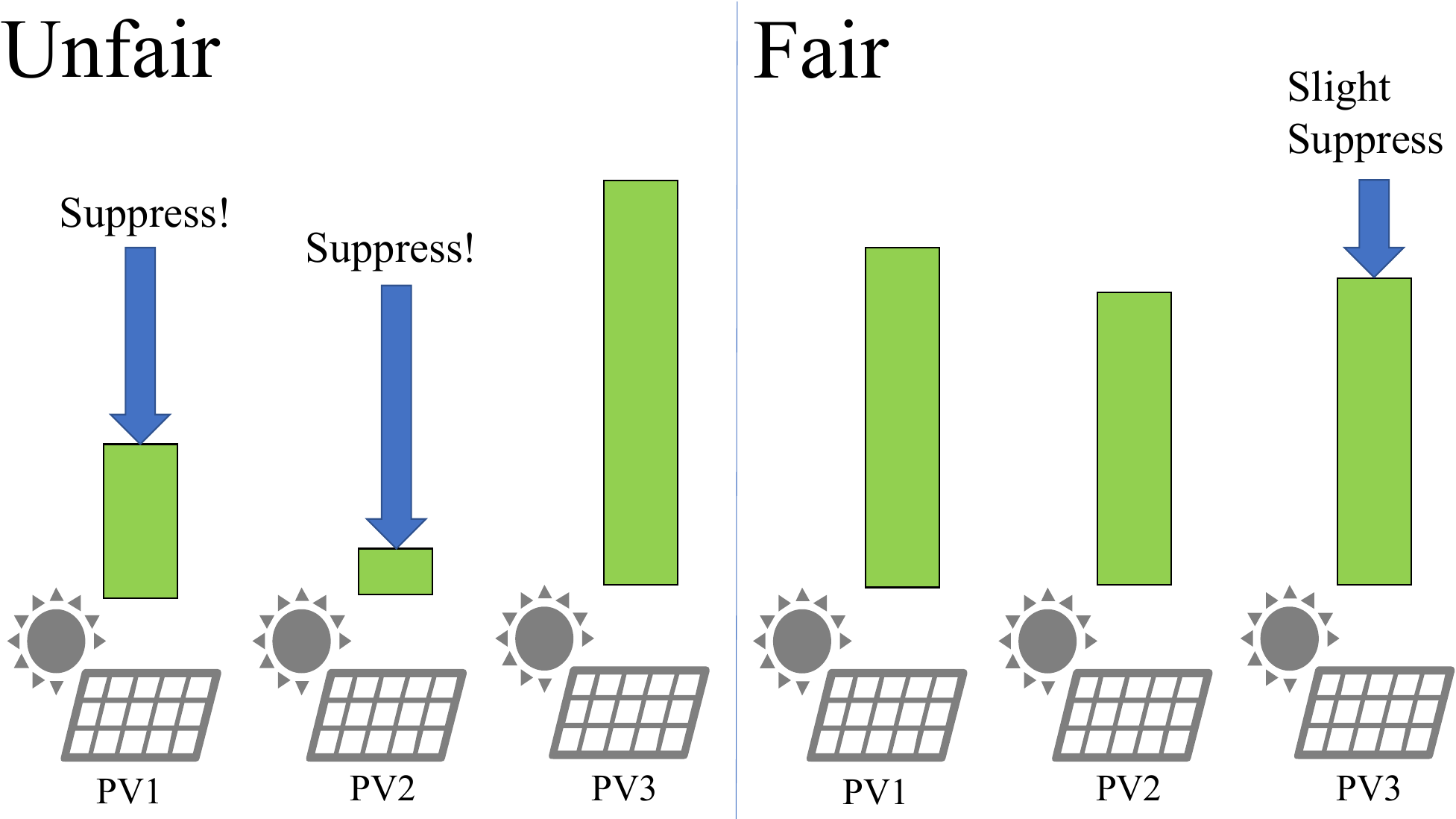}
\end{center}
\caption{Concept of Fairness among PV}\label{fig:Concept_of_Fairness_among_PVs}
\end{figure}

In the context of electricity deregulation, it is not uncommon for PV owners to be separate from the power company responsible for generating plans.
In such cases, the power company operator must be considerate of fairness among the owners, meaning that an equitable allocation of PVs is necessary as shown in Fig.\ref{fig:Concept_of_Fairness_among_PVs}.
This is because curtailment deprives generators of revenue opportunities, resulting in a loss of income for generator owners. Even if curtailment itself is unavoidable for system reliability, if generator owners perceive it as unfair, they will undoubtedly exit the market without hesitation.
However, our previous RE-RP algorithm does not incorporate fairness mechanisms and cannot always assign PV suppression fairly (just as it cannot always be unfair).
Now, we need to expand RE-RP to handle equal opportunity among PVs.
Some works related to fairness exist.
In reference~\cite{UC-equality}, the conventional UC problem is studied.
In reference~\cite{PVsupVoltVar} and~\cite{PVsupOPF}, 
the volt/var control and the optimal Power Flow (OPF) are respectively investigated.
However, in the robust optimization field, we have not seen any applications tackling fairness concept.
A more detailed discussion on comparison with the existing works~\cite{UC-equality}~\cite{PVsupVoltVar}~\cite{PVsupOPF} using fairness measures is given in Section\ref{subsec:Modeling of Fairness}.

This paper is organized as follows:
Section~\ref{sec:preliminaries} presents some mathematical preliminaries and introduces the deterministic model, which is in common with the RE-RP formulation.
In Section~\ref{sec:rerp}, after presenting the modeling of equal opportunity concept,
 we propose a new model \textit{RE-RP with Fairness~(RE-RP\textit{fair})}, which is an expansion from RE-RP.
Furthermore, we prove that Benders Decomposition algorithm can solve the new model.
In Section~\ref{sec:ExperimentalResults}, we see the effectiveness of RE-RP\textit{fair} through experiments.
Using artificial setting, we confirm that RE-RP\textit{fair} enhances fairness.
In Section~\ref{sec:conclusion} we give conclusion with future research directions.
\section{Preliminaries}\label{sec:preliminaries}
We first introduce the deterministic problem as same as RE-RP.
Then, to define RE-RP\textit{fair}, we add a penalty cost term to the objective function of the robust problem.
\subsection{The Deterministic Problem}
First, we suppose that the uncertain parameters are given 
as the load $d$, at the different nodes, and the output $z$ of the PVs for a deterministic problem
\footnote{
Definition of the symbols are in ~\ref{Appendix:Nomanclature}
}
.
Then, the problem can be modeled by the following mixed-integer linear optimization program:
\begin{flalign}
\min\limits_{x,u,v,r,p,q} & \sum\limits_{t=1}^T\sum\limits_{i=1}^{N_g}x^t_i F^t_i+u^t_i S^t_i+v^t_i G^t_i+C^t_i p^t_i+\sum\limits_{l=1}^{N_p} \Pi^t_lr^t_l\\
\mbox{subject to} \ &x^{t-1}_i-x^t_i +u^t_i  \ge 0,\qquad  \forall i \in \mathcal{N}_g,t\in \mathcal{T} \\
 &x^{t}_i-x^{t-1}_i +v^t_i  \ge 0,\qquad  \forall i \in \mathcal{N}_g,t\in \mathcal{T}  \\
&x^{t}_i-x^{t-1}_i \le x^{\tau}_i \qquad \forall \tau \in [t+1,\min\{t+\mbox{MinUp}_i \notag\\
& \hspace{30mm}-1,T\}],t\in[2,T]\\
& x^{t-1}_i-x^{t}_i \le 1-x^{\tau}_i \qquad \forall \tau \in [t+1,\min\{t \notag \\
&  \hspace{30mm} +\mbox{MinDw}_i-1,T\}], \notag \\
& \hspace{45mm} t\in[2,T]\\
& \sum\limits_{i=1}^{N_g}p^t_i +\sum\limits_{l=1}^{N_p}(1-r^t_l)\bar{z}^t_l=\sum\limits_{j=1}^{N_d}\bar{d}^t_j,  \forall t\in \mathcal{T}\\
&-RD^t_i \le p^t_i -p^{t-1}_i \le RU^t_i, \  \forall i \in \mathcal{N}_g, t\in \mathcal{T}\\
& p^{\min}_i \le p^t_i + q^t_{i}\le p^{\max}_i, \hspace{2mm} \forall i \in \mathcal{N}_g, t\in \mathcal{T} \\
 & p^{\min}_ix^t_i \le p^t_i \le p^{\max}_ix^t_i, \hspace{2mm} \forall i \in \mathcal{N}_g, t\in \mathcal{T} \\
&\sum\limits_{i \in \mathcal{N}_g}q^t_{i} \ge \bar{q}^t, \hspace{25mm} \forall t \in \mathcal{T} \\
&q^t_{i} \le \bar{q}^t_{i}, \hspace{21mm} \forall i \in \mathcal{N}_g, t \in \mathcal{T} \\
&x^t_i, u^t_i, v^t_i \in \{0,1\}, \hspace{8mm} \forall i \in \mathcal{N}_g, t\in \mathcal{T} \\
&r^t_l \in \{0,1\}, \hspace{16mm}  \forall l \in \mathcal{N}_p, t\in \mathcal{T} \\
 &p^t_i \ge 0, \hspace{22mm}  \forall i \in \mathcal{N}_g, t\in \mathcal{T}.
\label{eq:lp2}
\end{flalign}
The deterministic problem we consider is similar to the one introduced in~\cite{Bertsimas} and can be written, in a compact matrix formulation as:
\begin{flalign}
\begin{array}{llll}
\min\limits_{x,y} & \transpose{c} x + \transpose{b}y  \\
\mbox{subject to }& Fx \le f \\
& Hy \le G_1\bar{d}+G_2\bar{z}+h \\
& Ax+By \le g \\
& I_u y = \bar{d} -\bar{z} \\
& x \in \{0,1\}^n,\ y\ge 0.
\end{array}
\label{eq:lp}
\end{flalign}
The vector $x=(x_1,r)$ contains a subvector $x_1$ which is a vector of commitment related decisions including the on/off and start-up/shut-down of each generation unit for each time interval (as in \cite{Bertsimas}).
$r$ is a vector related to the usage of the solar panels at each time interval.
The vector $y$ is a vector of continuous variables, as in ~\cite{Bertsimas}, including the generation output, load consumption levels, resource reserve levels and power flows in the transmission network for each time interval.
 The term $\transpose{c}x$ in the objective function is related to the commitment cost while $\transpose{b}y$ is the dispatch cost over the planning horizon.
 The vector $\bar{d}$ contains the expected demand at the nodes at the different times and the vector $\bar{z}$ contains the expected output of a solar panel at the different nodes at the different times.

\section{RE-RP with Fairness}\label{sec:rerp}
In this section, we propose a new model, RE-RP with fairness (RE-RP\textit{fair}), that achieves both robust optimization and fair allocation among PVs.

\subsection{Modeling of Fairness}\label{subsec:Modeling of Fairness}
First, we propose a mathematical formulation of fairness. In previous papers, several fairness models have been adopted. For example, ~\cite{UC-equality} used the number of loss opportunities, and ~\cite{PVsupVoltVar,PVsupOPF} applied Jain's Fairness Index (JFI)~\cite{JFI} to the set of possible power output generated by the PVs.

Although there are various candidates for the fairness measure, we have chosen the Gini Index~\cite{Gini_Lorenz} because it has been widely used as a metric of income inequality in social science studies. Under deregulation, the PVs' output we focus on can be directly transferred to the income of the PVs' owner. Therefore, we apply the Gini Index to our research evaluation.

However, unfortunately, the Gini Index cannot be applied to RE-RP as it is because a sorting process is required in its calculation. It is not compatible with mathematical optimization techniques. Therefore, we introduce the sum of absolute deviations of each PV's output from its mean value as an alternative to the Gini Index. It has two good aspects for our purpose. First, it is known as one of the mathematical representations of diversity. If all PVs' power output is the same, the term is equal to zero. Second, the absolute value function can be reduced to linear terms, making it tractable in mathematical optimization theory. Therefore, we can solve RE-RP\textit{fair} in a similar way to how we would solve it with the Gini Index.

To summarize, we use the sum of absolute deviations in the optimization algorithm, and we also use the Gini Index in the evaluation of the algorithm.

To realize fairness allocation, we model the quantity of inequality among PVs as the sum of the absolute value of each PVs deviation from its mean value (\emph{e.g.}, overall a day).
We introduce new variables $s_l$ and $L^{1}(s)$, where 
$s_l$ models the total power of PV $z_l$, i.e.,
\begin{align}\label{eq:total power}
s_l = \sum_{t \in \tau}^{N_p} z_l^t,
\end{align}
and where $L^1(s)$
is defined by
\begin{align}
L^1(s)  = \sum\limits_{j = 1}^{N_p} \left| \frac{\sum\limits_{k = 1}^{N_p} s_k}{N_p} - s_j \right|.
\end{align}
To linearize the terms with absolute value, we need two new constraints.
\begin{align}\label{eq:L1const1}
& \frac{1}{N_p} \left( \sum\limits_{k = 1}^{N_p}\sum\limits_{t = 1}^{\mathcal{T}}((\bar{z}^t_k- \eta^t_k \hat{z}^t_k)(1-r^t_k)) \right)  \notag \\
 &- \sum_{t=1}^{\mathcal{T}} ((\bar{z}^t_j- \eta^t_j \hat{z}^t_j)(1-r^t_j)) =  \ a^+_j - a^-_j  \ \ \ (\forall j),
\end{align}
\begin{align}\label{eq:L1const2}
a^+_j , a^-_j \geq 0  \ \ \ (\forall j).
\end{align}
We can then rewrite $L^1(s)$ as $L^1(a^+, a^-)$,
\begin{align}
L^1(a^+, a^-) & = \sum\limits_{j = 1}^{N_p} ( a^+_j +  a^-_j \ )   \ \ \ (\forall j).
\end{align}

Next, we add $\chi L^{1}(a^{+}, a^{-})$
to the objective function of the deterministic problem as a penalty function. Here, $\chi $ is a positive parameter.
Adding (\ref{eq:L1const1}) and (\ref{eq:L1const2}) leads to the following objective function for RE-RP\textit{fair}:
\begin{align}
\min\limits_{x,u,v,r,p,q} & \sum\limits_{t=1}^T\sum\limits_{i=1}^{N_g}x^t_i F^t_i+u^t_i S^t_i+v^t_i G^t_i+C^t_i p^t_i 
 +\sum\limits_{l=1}^{N_p} \Pi^t_lr^t_l  \notag \\ & + \chi L^1(a^+, a^-).
\label{eq:newobj}
\end{align}
If $L^1(a^+, a^-)$ (deviation of PVs) increases, the objective function value will also increase.
Hence, since we are considering a minimization problem, RE-RP\textit{fair} will return an optimal solution with a small value for $L^1(a^+, a^-)$,
implying fairness among PVs production.

\subsection{Two-stage Adaptive Robust Formulation for RE-RP\textit{fair}}

\subsubsection{The Uncertainty Set}
The demand $d$ at the different nodes as well as the output $z$ of the solar panels are uncertain. In this paper, we model the uncertainty by assuming that the demand $d^t_j$, at time $t$ and node $j$, and the solar panel output $z^t_i$ at time $t$ and node $i$ can be written as 
\begin{align}
&d^t_j= \bar{d}^t_j+ \zeta^t_j \hat{d}^t_j \\
&z^t_j= (\bar{z}^t_j- \eta^t_j \hat{z}^t_j)(1-r^t_j)\ \ \mbox{ }
\end{align}
(notice here that $z^t_j$ is non-zero only if $r^t_j=0$).
Here, $\zeta^t_j \in [-1,1]$ and $\eta^t_j \in [-1,1]$ are variables that model the uncertainty respectively on $d$ and on $z$ and are assumed to belong to the uncertainty sets:
\begin{align}
\forall t,\ \zeta^t \in \mathcal{D}^t=\{\zeta^t\ |\ \sum_j |\zeta^t_j| \le \Delta^t,\ \zeta^t_j \in [-1,1]\}, \\
\forall t,\ \eta^t \in \mathcal{Z}^t=\{\eta^t\ |\ \sum_i |\eta^t_i| \le \Gamma^t,\ \eta^t_j \in [-1,1]\}.
\end{align}
We introduce $\mathcal{D}=\prod_t \mathcal{D}^t$ and $\mathcal{Z}=\prod_t \mathcal{Z}^t$.

\subsubsection{The Robust Formulation}
We now present a two-stage robust formulation of RE-RP\textit{fair}. We aim to find a vector $x$ that will minimize the cost 
$\transpose{c}x + \transpose{b}y(x,\zeta,\eta)$
in the ``worst possible'' scenario. Here $y(x,\zeta,\eta)$ is the optimal second-stage decision which depends both on $x$ and on the uncertainty variables $\zeta$ and $\eta$. Here ''worst'' means that for each $x$ we look for the scenario (i.e. the value of $\zeta$ and $\eta$) that will maximize the cost $\transpose{b}y(x,\zeta,\eta)$.
Further a linear penalty $-\chi L^1(a^+, a^-)$ (i.e., convex) is added to the objective.
Hence the two-stage adaptive robust problem can be formulated as 
\begin{align}
&\min_x (\transpose{c}x  \notag \\
&+\max\limits_{\zeta,\eta \in \mathcal{D}\times \mathcal{Z},a^+,a^-} (-\chi L^1(a^+, a^-)+ 
&\min\limits_{y \in \Omega(x,\zeta,\eta)} \transpose{b}y ) ) \\
&\mbox{subject to } Fx \le f, \\
&x=(x_1,r) \in \{0,1\}^n, 
\end{align}
where
$$\Omega(x,\zeta,\eta)=\{y\ |\ Hy \le G_1 d+G_2{z}+h,\ I_u y = \bar{d} -\bar{z},\ Ax+By \le g  \}. $$
Let us denote by $S(x,\zeta,\eta)$ the optimal value of the following linear optimization problem:
$$ S(x,\zeta,\eta)=\min\limits_{y \in \Omega(x,\zeta,\eta)}\transpose{b}y,  $$
which corresponds to the economic dispatch cost when $(x,\zeta,\eta)$ is fixed. Since we maximize $S(x,\zeta,\eta)$ over $(\zeta,\eta) \in \mathcal{D}\times \mathcal{Z}$ the worst case will occur when both $\zeta,\eta \ge0$. Hence we can simplify the uncertainty sets to
\begin{align}
\forall t,\ \zeta^t \in \mathcal{D}^t=\{\zeta^t\ |\ \sum_j \zeta^t_j \le \Delta^t,\ \zeta^t_j \in [0,1]\}, \\
\forall t,\ \eta^t \in \mathcal{Z}^t=\{\eta^t\ |\ \sum_i \eta^t_i \le \Gamma^t,\ \eta^t_j \in [0,1]\}. 
\end{align}

By duality we  have
\begin{align}
 S(x,\zeta,\eta)= \max\limits_{\phi,\lambda,\mu}\  & \transpose{\lambda}(Ax-g)-\transpose{\phi}(G_1 (\bar{d}+ \zeta \circ \hat{d} ) \notag \\
 &+ G_2{(\bar{z}+ \eta \circ \hat{z} )\circ r}+h)+ \transpose{\mu}((\bar{d}+ \zeta \circ \hat{d} ) \notag \\
 &- (\bar{z}+ \eta \circ \hat{z} )\circ (1-r)) \\
 \mbox{subject to}& -\transpose{\lambda}B - \transpose{\phi}H + \transpose{\mu}I_u = \transpose{b} \\
 &\phi , \lambda\ge 0,
\end{align}
where $a \circ b$ denotes the vector having the same size as $a$ and $b$ and whose $i$th component is $a_i b_i$. \\

Let us denote by $R(x)$ the optimal value of the recourse problem (i.e. $R(x)=\max\limits_{\zeta,\eta,a^+,a^-}\left(-\chi L^1(a^+, a^-)+S(x,\zeta,\eta)\right)$), we have
\begin{align}
R(x)= \max\limits_{\phi,\lambda,\mu,\zeta,\eta,a^+,a^-}\  & -\chi L^1(a^+, a^-)+\transpose{\lambda}(Ax-g) \notag \\ 
 &-\transpose{\phi}(G_1 (\bar{d}+ \zeta \circ \hat{d} ) \notag \\ 
 &+G_2{(\bar{z}+ \eta \circ \hat{z} )\circ r}+h) \notag \\ 
 &+\transpose{\mu}((\bar{d}+ \zeta \circ \hat{d} ) \notag \\
 & -(\bar{z}+ \eta \circ \hat{z} )\circ (1-r)) \\
\mbox{subject to}& -\transpose{\lambda}B - \transpose{\phi}H + \transpose{\mu}I_u = \transpose{b} \\
&\phi , \lambda\ge 0,  (\zeta,\eta) \in \mathcal{D}\times \mathcal{Z}. \label{eq:RXlast}
\end{align}
The goal now is to show that we can reformulate the objective function of $R(x)$ into a convex function.
$-\chi L^1(a^+, a^-)$ is a convex function by definition of $L^1(\cdot)$, however the remaining terms, i.e., $\transpose{\lambda}(Ax-g)-\transpose{\phi}(G_1 (\bar{d}+ \zeta \circ \hat{d} )+G_2{(\bar{z}+ \eta \circ \hat{z} )\circ r}+h)+ \transpose{\mu}((\bar{d}+ \zeta \circ \hat{d} )-(\bar{z}+ \eta \circ \hat{z} )\circ (1-r))$ are problematic as the bi-linear terms $\transpose{\phi}(\zeta \circ \hat{d} ),\ \ \transpose{\phi}(\eta \circ \hat{z}\circ r) $,  $\transpose{\mu} (\zeta \circ \hat{d} )$ or $\transpose{\mu}(\zeta \circ \hat{d} \circ(1-r))$ are not convex. Nonetheless we will show that these terms can  be linearized by noticing that:

\begin{itemize}
	\item For any $x=(x_1,r)$ there exists  an optimal solution $(\zeta^*,\eta^*)$, which is a binary vector.
	Indeed, if we fix $\phi,\lambda,\mu$, in the objective function of $R(x)$, the function is convex. Indeed $-\chi L^1(a^+, a^-)$ is convex, and if $\mu$ and $\phi$ are fixed, the problematic bi-linear terms described above become linear terms (for example $\transpose{\phi}(\zeta \circ \hat{d} )$ is a linear term when $\phi$ is fixed).
 Since $(\zeta^*,\eta^*)$ are decoupled from the other variable we see that a maximum value in $(\zeta^*,\eta^*)$ is attained at an extreme point of $\mathcal{D}\times \mathcal{Z}$ (by convexity of the function). Hence we can assume that $(\zeta^*,\eta^*)$ is a binary vector.\\
	\item By the point above, the problematic bi-linear terms, $\transpose{\mu} (\zeta \circ \hat{d} )$ and $\transpose{\phi}(\zeta \circ \hat{d} )$  are a product of continuous (dual) variable ($\phi$ or $\mu$) by a binary variable ($\zeta$ or $\eta$),  and hence can be rewritten as the product $\kappa \beta$ of a positive continuous variable $\kappa$ by binary variable $\beta$. This product can be linearized by introducing a new variable $v= \kappa \beta$ that satisfies the following linear constraints:
	\begin{align*}
	& v \le \kappa_m \beta \\
	& v \le \kappa \\
	& v \ge \kappa - \kappa_m(1-\beta) \\
	& v \ge 0, \beta \in \{0,1\}, \kappa \ge 0,
	\end{align*}
	where $\kappa_m$ is an upper bound on $\kappa$.
\end{itemize}

Now that the objective function of $R(x)$ has been convexified, we now that for every $x$, there exists an optimal solution of $R(x)$ at an extreme point of the feasibility region (that is because we are maximizing a convex function). More precisely, let $\Xi$ denote the set of constraints of the dual problem $S(x,\zeta,\eta)$ and let $\Lambda$ denote the set of constraints induced by the linearizations. Furthermore, let us denote by $F(x)$ the linear part of the objective function of $R(x)$ and by $\tilde{y}$ the aggregated variables (the dual variables, the uncertainty variable $\zeta,\eta$ and the new variables induced by the linearization). We rewrite $R(x)$
\footnote{Overall of $R(x)$ is in \ref{Appendix:Recourse Problem}}
as:

\begin{flalign}
\max\limits_{\eta,\zeta,\tilde{y},a^+,a^-} & -\chi L^1(a^+, a^-)+ F(x)^\top \tilde{y} \\
\mbox{subject to }& \tilde{y}\in \Xi\cap \Lambda \\
& (\zeta,\eta) \in \mathcal{D}\times \mathcal{Z}\\
& \zeta, \eta \in \{0,1\}.
\end{flalign}

Let $\mathcal{C}$ denote the convex hull of the extreme points in the feasible region of $ RP(x)$ (notice that the feasible region is independent of the value of $x$)
and let $\mathcal{S}=\{(\tilde{y}^s,\zeta^s,\eta^s),s=1,\dots,S\}$ denote the set of extreme points of $\mathcal{C}$. Since, for any $x$, the objective function of $R(x)$ is convex, there exists always a point of $\mathcal{S}$ that is an optimal solution of $R(x)$.
The robust problem $RP$
\footnote{Overall of $RP$ is in \ref{Appendix:Robust Problem}}
can be rewritten as:
\begin{flalign}
\min\limits_{x,w} & ~c^\top x+w\\
\mbox{subject to }& Fx \le f\\
& x\in \{0,1\} \\
& w \ge -\chi L^1(a^+s, a^-s)+ F(x)^\top \tilde{y}^s, \hspace{3mm} \forall \tilde{y}^s \in \mathcal{S}. 
\end{flalign}
The problem can now be solved using Benders decomposition algorithm Algorithm~\ref{ALG:BendersALG}  as in~\cite{Bertsimas}. 
\subsection{Algorithm and Proof}
\begin{figure}[tb]
\begin{center}
 \includegraphics[width=9.0cm]{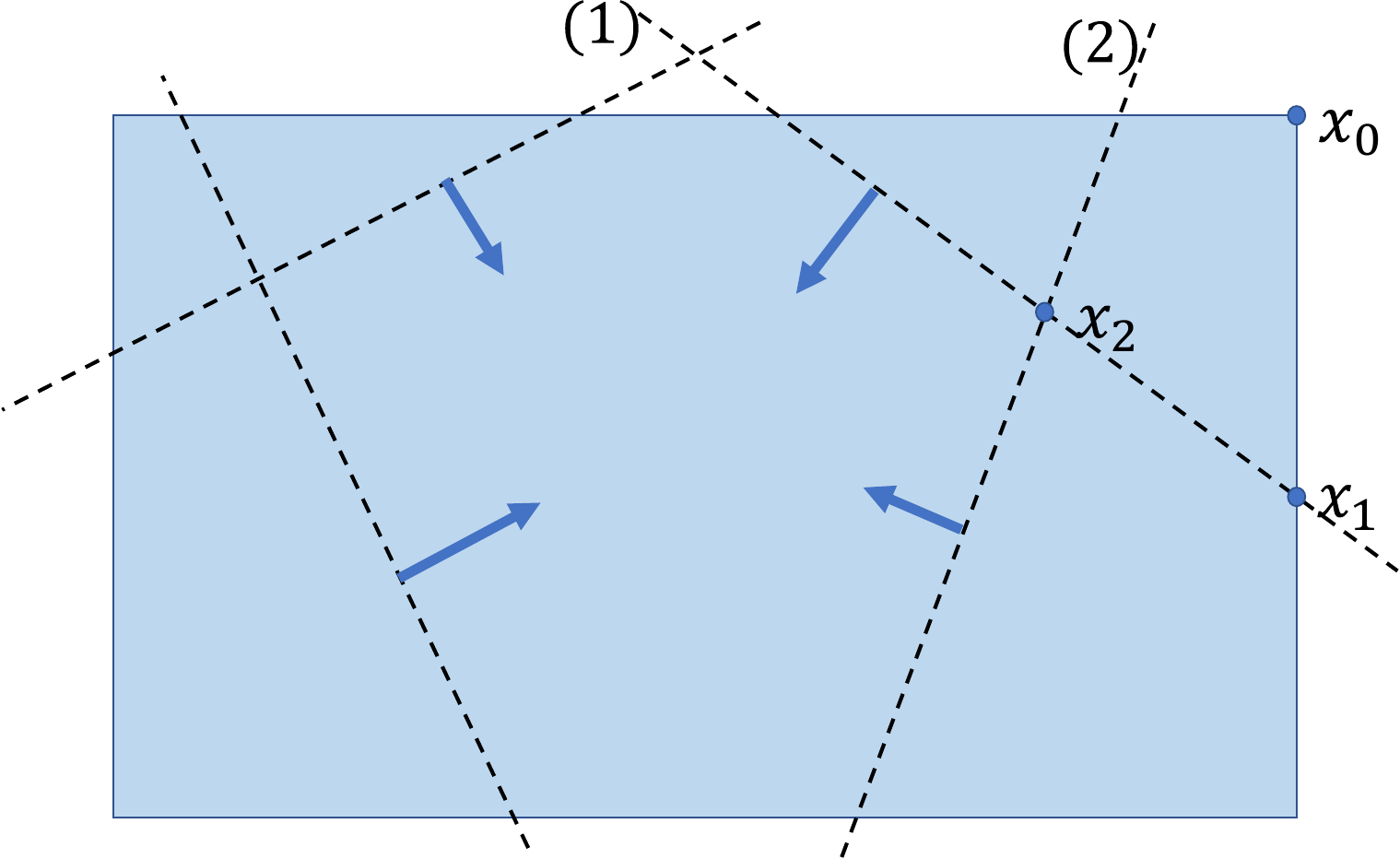}
\end{center}
\caption{2-iterations of Benders algorithm: the blue region represents the constraint $Fx \leq f$ and the dashedlines	represent the inequalities added at each iteration.}\label{fig:BendersImage}
\end{figure}

Now, we show Algorithm~\ref{ALG:BendersALG} as Benders Decomposition Algorithm for RE-RP\textit{fair}, and prove its correctness.
\begin{algorithm}[t]
\caption{Benders decomposition algorithm for RE-RP\textit{fair}}
 
 \begin{algorithmic}[1]\label{ALG:BendersALG} 
  \STATE $U_0=\infty$, $k=1$ (initialization)

  \STATE Solve the following optimization problem, denote its optimal solution by $(x^k,w^k)$ and its optimal value by $L_k$.
\begin{align}
RP : \min_{x,w} & \left(\transpose{c}x + w\right) \notag \\
\mbox{subject to }& Fx \le f \notag\\
& x \in \{0,1\}^n, w \ge 0\notag
\end{align}

  \STATE Solve the recourse problem, $R(x^k)$, Eq.(~\ref{eq:RXlast})~when $x=x^k$ and denote its optimal value by $v(R(x^k))$.
  
  \STATE $U_k = \transpose{c}x^k + v(R(x^k))$
  
  \STATE If $(\inf_{j \le k}{U_j}-L_k)/(\inf_{j \le k}{U_j})<\epsilon$, solve the linear optimization as shown below
\begin{align}
S(x,\zeta,\eta):\min\limits_{y \in \Omega(x,\zeta,\eta)}\transpose{b}y,
\end{align}
and terminate the algorithm.

Otherwise, set $k=k+1$, return to  step~2 and add the following constraint to the  problem RP.
\begin{align}
\label{eq:wconst}
w \ge  & -\chi L^1(a^+k, a^-k) + \transpose{\lambda^k}(Ax-g)-\transpose{\phi^k}h \notag \\ &+ \transpose{\mu^k}((\bar{d}+ \zeta^k \circ \hat{d} ) -(\bar{z}- \eta^k \circ \hat{z} )\circ (1-r)), 
\end{align}
where $(\phi^k,\lambda^k,\mu^k,\zeta^k,\eta^k,a^+k,a^-k)$ is an optimal solution of the recourse problem when $x=x^k$.
\end{algorithmic}
\end{algorithm}
\begin{theorem}\label{LkMonotoInc2}
Algorithm 1 obtains $v(RP^*)$ after a finite number of steps.
\end{theorem}
\textit{Sketch of the proof:}\label{L_k converge}
We recall the equivalent linear reformulation of the robust problem.
\begin{align}
RP : \min_{x,w} & \left(\transpose{c}x + w\right)  \\
\mbox{subject to }  Fx &\le f \\
w &\ge	-\chi L^1(a^+k, a^-k) + \transpose{\lambda^k}(Ax-g) \notag \\
 &- \transpose{\phi^k}h+ \transpose{\mu^k}((\bar{d}+ \zeta^k \circ \hat{d} ) \notag\\
 &- (\bar{z}- \eta^k \circ \hat{z} )\circ (1-r))	\ \ \forall k \in \Psi (*)  \\
 x &\in \{0,1\}^n, w \ge 0. \notag
\end{align}
Here the set $\Psi$ indexes all the extreme points of the convex hull of the feasible region of the recourse problem $R(x)$.
The correctness of Algorithm 1 is proved in reference~\cite{roboptrecourse}, however for the sake of completeness of the paper we recall here the main ideas of the proof using Fig~\ref{fig:BendersImage} where the inequalities $(*)$ are represented by the dashed-lines.
Initially no inequality $(*)$ is known and the minimum is computed over the constraints
$Fx \leq f$ (the region in blue in the Fig.\ref{fig:BendersImage}).
At each step of the algorithm a new inequality (corresponding to the dashed-line in the figure) is added.
Since each inequality corresponds to an extreme point of the convex-hull of the feasibility region of $R(x)$, the number of these inequalities is finite.
Therefore after a finite number of iterations all the inequalities of the problem have been added and we find therefore the solution of RP
\footnote{If at the $k$-iteration the algorithm returns a solution $x_k$ that satisfies all the inequalities $(*)$ (even the inequalities that remain to be added yet) then $x_k$ is an optimal solution of RP.}
.
\QED
\section{Experimental Results}\label{sec:ExperimentalResults}
In this section, we explain the effectiveness of  RE-RP\textit{fair}, through numerical experiments.
\subsection{Settings}
For the experiments, we built an artifitial model of thermal generators, electricity demand, and PV outputs. We assume this model a power system on an isolated island. However, we refer to Miyako Island~\cite{MIYAKOJIMA} as the peak of the demand and the sum of the PV and thermal generator output. First, we set the basic parameters to $N_g=3$, $N_d=10$, $N_p=3$, and $T=24$.

Table~\ref{table:Thermal Generator Parameters} shows the parameters of the thermal generators.
\begin{table}[bh]
  \begin{center}
 \caption{Thermal Generator Parameters}
 \label{table:Thermal Generator Parameters}
   \renewcommand{\arraystretch}{1.2}
  \tabcolsep=11pt 
 \begin{tabular}{c|c|c|c} \hline
      Constant & G1 & G2 & G3(Slack)  \\ \hline \hline
      F[kJPY] & 0.0 & 0.0 & 10,000 \\
      S[kJPY] & 25 & 6.3 & 10,000 \\
      G[kJPY] & 10 & 2.5 & 0 \\
      C[kJPY/MW] & 13.4 & 26.0 & 10,000 \\
      RD[MW/h] & 60.0 & 15 & 100,000 \\
      RU[MW/h] & 60.0 & 15.0 & 100,000 \\
      MinUp[h] & 6.0 & 3.0 & 0.0 \\
      MinDw[h] & 7.0 & 4.0 & 0.0 \\
      $P_{max}$[MW] & 60.5 & 10.5 & 100,000 \\
      $P_{min}$[MW] & 6.5 & 1.5 & 0 \\
      $\bar{q}^t_i$[MW] & 45.4 & 11.3 & 100,000 \\ \hline
 \end{tabular}
  \end{center}
\end{table}
G1 and G2 correspond to real generators, with G1 serving as the base power source due to its relatively low cost. G2 can be used during peak times or when G1 is out of service. On the other hand, G3 is a slack (virtual) unit with more parameters than the others to prevent infeasible solutions during iterations.

\begin{figure}[tb]
\begin{center}
 \includegraphics[width=8.0cm]{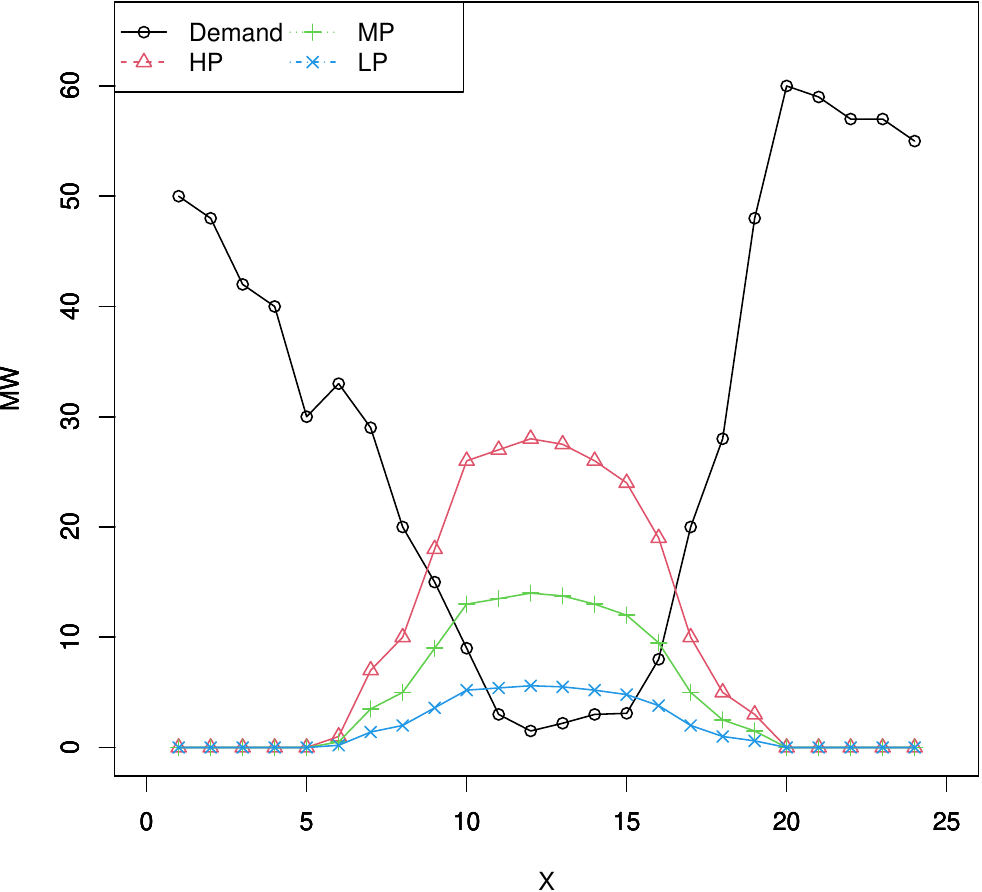}
\end{center}
\caption{Electricity Demand and PV Output Patterns}\label{fig:Electricity Demand and PV output Patterns}
\end{figure}

Next, we prepared the electricity demand and PV outputs. In a real situation, they would come from external forecast systems.
Fig.\ref{fig:Electricity Demand and PV output Patterns} shows the demand and PV patterns which are used in later experiments.
We use the only one demand pattern in all experiment in this paper, commonly.
It is sum of all individual demand $\bar{d}$ in each slot.
We allocate individual demand by distributing the demand pattern, randomly.
Furthermore, we introduce three PV Patterns.
HP means high-power output corresponding to sunny day.
On the other hand, LP represents low-power output as rainy day.
MP is middle-power representing medium case.
We also make individual PV output $\bar{z}$ likewise above demand case.
Uncertainty parameters $\hat{d}^{t}_{i},\hat{z}^{t}_{i}$ are set as below,
\begin{align}
&\hat{d}^t_j=  \ \ 0.2 * \bar{d}^t_j,  \\
&\hat{z}^t_j=  \ \ 0.2 * \bar{z}^t_j.
\end{align}
We have fixed the coefficient 0.2, though it can be tuned depending on the application. 

Finally, we fix the PV shutdown cost as $\Pi^{t}_{l} = 11.0 * \bar{z}^{t}_{l}$. The value of 11.0 Yen/kWh is referenced from the Japan Feed-in Tariff (FIT) rates in 2021~\cite{JAPAN-FIT}.

\subsection{Computational Parameters}
We use the Gurobi optimizer as subroutines of our Benders Decomposition. We use version 9.1.0 of Gurobi, which runs on a server with an Intel(R) Xeon Gold 6142 2.6GHz CPU and 64GB of RAM. The MIPGAP parameter is set to the default value of $10^{-4}$. In addition, we set the Benders Decomposition parameters such that $\epsilon$ to $0.001$ and the iteration limit is set to $30$.
\subsection{Evaluation Metric}
To evaluate fairness quantitatively, we introduce the metric called Gini Index~\cite{Gini_Lorenz}. The Gini Index is known as a summary statistic of inequality and is commonly used in social science studies \cite{GiniIndex}. It is calculated for a dataset such as income or property, and it ranges from 0 to 1 continuously. A value of 0 indicates a perfect equality situation (where everyone has the same amount of money), while a value of 1 indicates a situation where only one person has all the money.

In our study, we have focused on fairness over the total power of PV (See Eq.(\ref{eq:total power})).
To calculate the Gini Index, we 
sort $s_i$ by increasing order. The Gini Index of $s$  is given bellow:
\begin{align}
\sum^{N_p}_{i=1} \left( \frac{i}{N_p} - \frac{\sum^{i}_{j=1}s_j}{\sum^{N_p}_{j=k}s_k} \right).
\end{align}
\subsection{Preliminary Experiment}

\begin{figure}[tb]
\begin{center}
 \includegraphics[width=8.0cm]{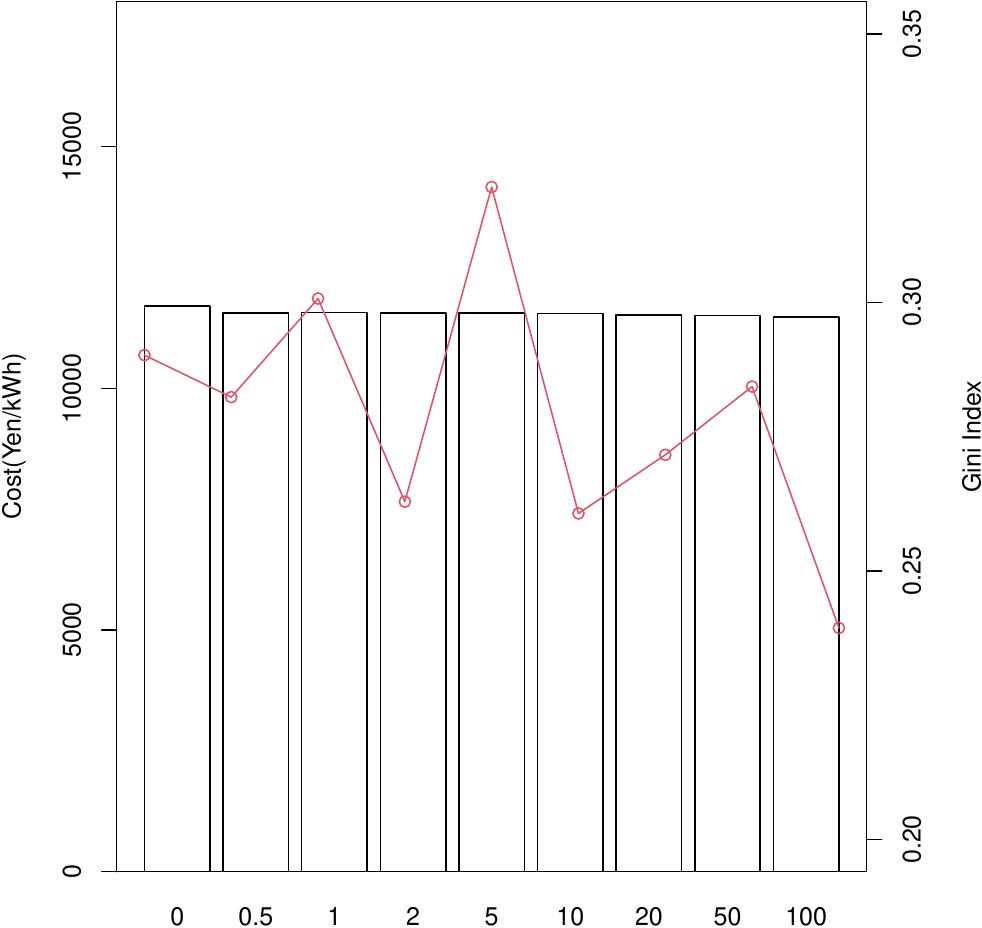}
\end{center}
\caption{Effects Of $\chi$ Parameter}\label{fig:Effects Of chi Parameter}
\end{figure}

The hyper-parameter $\chi$ plays a central role in the RE-RP$_\textit{fair}$ formulation. 
To evaluate the effects of $\chi$, we use the demand and HP case in Fig.\ref{fig:Electricity Demand and PV output Patterns}, showing
the results in Fig.\ref{fig:Effects Of chi Parameter}.
The $X$-axis in Fig.\ref{fig:Effects Of chi Parameter} 
represents the value of $\chi$ and the left $Y$-axis and the black bar plots show the cost including fuel cost and penalty for PV shutdown.
The right $Y$-axis is the value of Gini index.
We can observe the expected decreasing trend of the Gini Index. Additionally, we notice a similar trend in the cost. This contradicts our prediction, which is unfortunate. However, the fluctuations in the cost value are minimal compared to its absolute value. Therefore, we conclude that the negative impact of a larger $\chi$ is insignificant. Consequently, we set $\chi = 100.0$ for the main experiment.

\subsection{Main Results}
\begin{figure}[t]
\begin{center}
 \includegraphics[width=8.0cm]{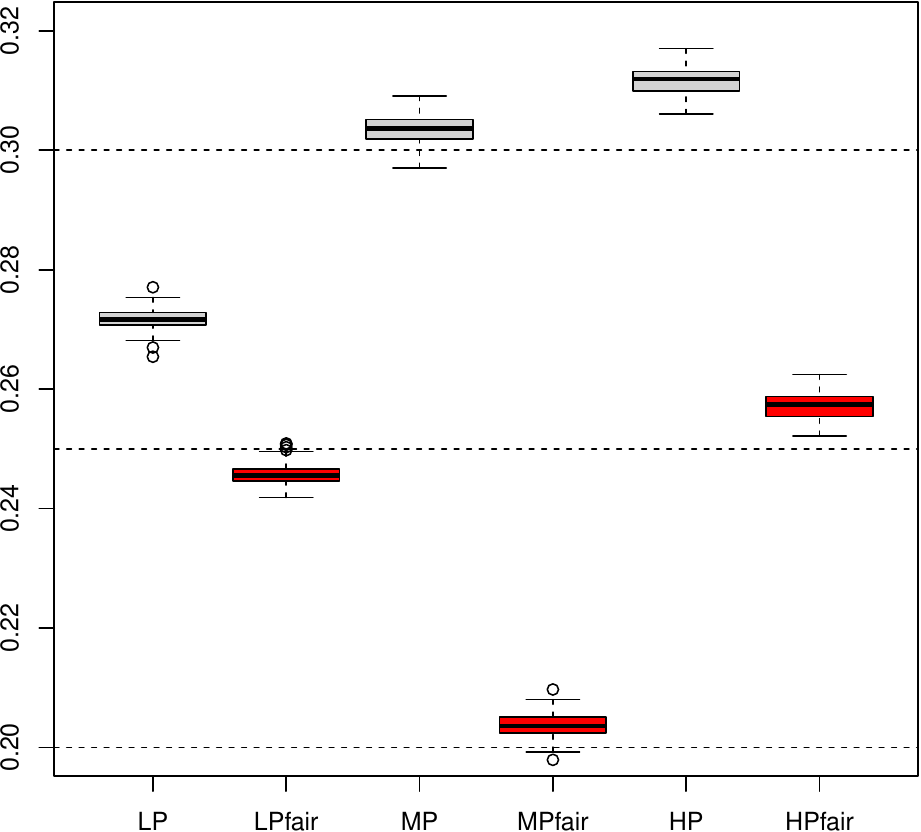}
\end{center}
\caption{Boxplot of the Gini Indices}\label{fig:Boxplot of the Gini Indices}
\end{figure}

To clarify the efficiency of RE-RP$_\textit{fair}$, we introduce six test cases composed by three pairs in Table~\ref{table:Teat Cases}.
In each pair, we check the Gini Index discrepancy between RE-RP$_\textit{fair}$ and existing RE-RP.
The Name column in Table~\ref{table:Teat Cases} represents PV patterns and optimization model.
For example, No.1 LP adopts PV pattern LP and RE-RP.
No.2 LP$_\textit{fair}$ takes same PV pattern and RE-RP\textit{fair}.
In MP, LP cases are likewise.
\begin{center}
\begin{table}[htb]
 \centering
 \caption{Test Cases}
 \label{table:Teat Cases}
 \scalebox{1.1}{
 \renewcommand{\arraystretch}{1.2}
 \tabcolsep=11pt
 \begin{tabular}{clcr}\hline
  \normalsize No & \normalsize Name & \normalsize PV & \normalsize $\Gamma$  \\ \hline
  1 & LP & Pattern LP & 0  \\ 
  2 & LP$_\textit{fair}$  & Pattern LP & 100  \\ \hline
  3 & MP & Pattern  MP& 0  \\ 
  4 & MP$_\textit{fair}$  & Pattern MP &100  \\ \hline
  5 & HP & Pattern HP & 0 \\ 
  6 & HP$_\textit{fair}$  & Pattern HP  & 100 \\ \hline
 \end{tabular}
 }
\end{table}
\end{center}
Next, we show how to generate data set in each case.
To show statistical validity, we simulate the forecast error by using random numbers.
In a real case,
our optimization algorithm would be used a day before. However,
in practice real PV output is determined on the day (i.e., operation day).
We define real PV output $z^{t}_{j,R}$, as 
\begin{align}\label{eq:z in real}
z^{t}_{j,R} \ = \ \bar{z}^{t}_{j} + N(\bar{z}^{t}_{j}, \frac{\hat{z}^{t}_{j}}{3}).
\end{align}
Then, we can calculate the total power $s_{j,R}$ using $z^{t}_{j,R}$, in the same way of (\ref{eq:total power}).
\begin{algorithm}
 \caption{Data Set Generation}
 \label{ALG:Data Set Generation}
 \begin{algorithmic}[1]
  \REQUIRE PV pattern and $\chi$
  \ENSURE  $M$ Set of Gini Indice
  \STATE Run Optimization
  \FOR{$k=1$ to $M$}
  \STATE Calculate $z^{t,k}_{j,R} (\forall t,j)$.
  \STATE Calculate $s^{k}_{j,R} (\forall j)$.
  \STATE Calculate Gini Index.
  \ENDFOR
 \end{algorithmic}
\end{algorithm}
Now, we introduce the main results.
Table~\ref{table:Mean and Standard Deviation of Gini Indices} shows the means and the standard deviations of the test cases.
In every pairs, RE-RP\textit{fair} is superior to RE-RP, as we expected.
To strengthen the conclusion, we explain the superiority in the statistical viewpoint by using hypothesis testing~\cite{StatisticalHypothesis}.
\begin{center}
\begin{table}[htb]
\centering
\caption{Mean and Standard Deviation of Gini Indices}
\label{table:Mean and Standard Deviation of Gini Indices}
\scalebox{1.1}{
\renewcommand{\arraystretch}{1.2}
 \tabcolsep=3pt
 \begin{tabular}{crrrrrr}\hline
  Case & LP & LP\textit{fair} & MP & MP\textit{fair} & HP & HP\textit{fair}  \\ \hline
Mean & 0.27 & 0.25 & 0.30 & 0.20 & 0.31 & 0.26  \\ 
Std & 0.0018 & 0.0019 & 0.0022 & 0.0019 & 0.0022 & 0.0021 \\ \hline
 \end{tabular}
}
\end{table}
\end{center}
First, we apply Shapiro-Wilk Test to check normality.
Table~\ref{table:Results of Shapiro-Wilk Test} shows w-value and corresponding p-value.
We can not reject null hypothesis; the distribution is normal, under significance level is 0.05, in all cases.
\begin{center}
\begin{table}[htb]
\centering
\caption{Results of Shapiro-Wilk Test}
\label{table:Results of Shapiro-Wilk Test}
\scalebox{1.1}{
 \renewcommand{\arraystretch}{1.4}
 \tabcolsep=3pt
 \begin{tabular}{crrrrrr}\hline
  Case & LP & LP\textit{fair} & MP & MP\textit{fair} & HP & HP\textit{fair}  \\ \hline
w-value & 0.97 & 0.99 & 0.99 & 0.99 & 0.99 & 0.99  \\ 
p-value & 0.052 & 0.37 & 0.552 & 0.734 & 0.721 & 0.969  \\ \hline
 \end{tabular}
}
\end{table}
\end{center}
Next, we check the discrepancy of the variances between pairs.
Table~\ref{table:Results of F Test} shows w-value and corresponding p-value.
We can not reject null hypothesis; each pair has same variance, under significance level is 0.05, in all pairs.
\begin{center}
\begin{table}[htb]
\centering
\caption{Results of F-Test}
\label{table:Results of F Test}
\scalebox{1.1}{
 \renewcommand{\arraystretch}{1.4}
 \tabcolsep=3pt
 \begin{tabular}{crrrrrr}\hline
  Case & LP vs  LP\textit{fair} & MP vs MP\textit{fair} & HP vs HP\textit{fair}  \\ \hline
F-value & 1.19 & 0.72 & 1.31  \\ 
p-value & 0.387 & 0.115 & 0.169 \\ \hline
 \end{tabular}
}
\end{table}
\end{center}
Finally, based on the results of Tables~\ref{table:Results of Shapiro-Wilk Test} and \ref{table:Results of F Test},
we apply Student's t-Test to check differences between the mean of each pair (Fig.\ref{table:Results of Student's t-Test}).
We can reject null hypothesis; each pair has the same mean, in all pairs.
Furthermore, we can see the differences in the boxplot shown in Fig.\ref{fig:Boxplot of the Gini Indices}.
Hence, RE-RP\textit{fair} can reduce the value of Gini index;
It enhances fairness allocation among PVs.
\begin{center}
\begin{table}[htb]
\centering
\caption{Results of Student's t-Test}
\label{table:Results of Student's t-Test}
\scalebox{1.1}{
\renewcommand{\arraystretch}{1.4}
 \tabcolsep=3pt
 \begin{tabular}{crrrrrr}\hline
  Case & LP vs  LPw/EO & MP vs MPw/EO & HP vs HPw/EO  \\ \hline
t-value & -96.02 & 51.35 & 110.82  \\ 
p-value & $<$ 0.001 & $<$ 0.001 & $<$ 0.001 \\ \hline
 \end{tabular}
}
\end{table}
\end{center}
\clearpage

\section{Conclusion}\label{sec:conclusion}
In this paper, we propose a new robust optimization model, RE-RP\textit{fair} to enhance fairness among PVs outputs.
RE-RP\textit{fair} has a penalty term composed of the sum of the absolute values of deviations of PV outputs.
It looks more complex than existing RE-RP model.
But we gave the proof that the same Benders Decomposition algorithm can solve RE-RP\textit{fair}.
Furthermore, we confirmed that RE-RP\textit{fair} can make PV outputs allocation more fair
for small island grid; RE-RP\textit{fair} reduces Gini indices in all cases.

As in the above-mentioned paper, we can not help using artificial data in the experiments.
Due to privacy and security concerns, it is hard to obtain each electricity demand and PVs outputs data.
While our experiments show the effectiveness of RE-RP\textit{fair}, we may need more realistic data to build more reliable model and implementation to ensure its reliability in real operation.

In the future, we would like to improve and tune our model to deal with pragmatic situations. We plan to tackle large-scale realistic settings to reveal the relationship between fairness and total cost.
\section{Acknowledgement}
The authors thank Bruno F. Lourenço for his helpful advice.
We also thankful to Taichi Isogai and Yuki Sato for supporting this research project.
\clearpage
\bibliographystyle{IEEEtran}
\bibliography{toyoshima}
\appendix
\clearpage
\section{Nomanclature}\label{Appendix:Nomanclature}
The decision variables in this paper are:
\begin{itemize}
	\item $x^t_i\in \{0,1\}$: if generator $i$ is on at time $t$, $x^t_i=1$ otherwise
	      $x^t_i=0$.
	\item $u^t_i\in \{0,1\}$: if generator $i$ is turned on at time $t$, $u^t_i=1$ otherwise $u^t_i=0$.
	\item $v^t_i\in \{0,1\}$: if generator $i$ is turned down at time $t$, $v^t_i=1$ otherwise $v^t_i=0$.
	\item $r^t_l \in \{0,1\}$: if PV $l$ is turned down at time $t$, $r^t_l=1$ otherwise $r^t_l=0$.
	\item $p^t_i \ge 0$: production of generator $i$ at time $t$.
	\item $q^t_{i} \ge 0$: reserve of generator $i$, time $t$
\end{itemize}
The parameters and sets associated to the problem are as follows:
\begin{itemize}
	\item $N_g,N_p,N_d,T$: the number of generators, photovoltaics (PVs), loads and time periods (in hours).
	\item $\mathcal{N}_g,\mathcal{N}_p,\mathcal{N}_d,\mathcal{T}$: the corresponding sets of generators, PV, loads  and time periods (in hours).
	\item $S^t_i, G^t_i, F^t_i$: start-up, shut-down, and no-load costs of generator $i$ at time $t$.
	\item $C^t_i, \Pi^t_l$: respectively the cost per MWh of generator $i$ at time $t$ 
          and the cost per MWh related to the usage of the PV $l$ at time $t$. 
	\item $p^{\max}_i,p^{\min}_i$: maximum and minimum production levels of generator $i$.
	\item $RU^t_i,RD^t_i$: ramp-up and ramp-down rates of generator $i$ at time $t$.
	\item $\mbox{MinUp}_i, \mbox{MinDw}_i$: minimum-up and minimum-down times of generator $i$.
	\item $\bar{d}^t_j,\bar{z}^t_l$: expected load at node $j$ at time $t$ and expected PV output at node $l$ at time $t$.
	\item $\bar{q}^t_{i}$: reserve capacity of generator $i$, time $t$.
	\item $\bar{q}^t$: reserve capacity requirement of power system at time $t$.
	\item $x^0_i\in \{0,1\}$: initial ``on/off'' state of generator $i$.
	\item $p^0_i \ge 0$: initial production of generator $i$.
\end{itemize} 
\section{Recourse Problem}\label{Appendix:Recourse Problem}
\begin{flalign}
\max  &- \Gamma \left( \sum\limits_{j = 1}^{N_p} ( a^+_j + a^-_j \ ) \right) \notag \\
& + \sum\limits_{t\in \mathcal{T}} (\sum\limits_{j=1}^{N_d}\bar{d}^t_j - \sum\limits_{j=1}^{N_p}(1-r^t_j)\bar{z}^t_j)(\alpha^+_t-\alpha^-_t) \notag \\
&+\sum\limits_{t\in \mathcal{T}}(\sum\limits_{j=1}^{N_d}\hat{d}^t_j(m^1_{j,t}-m^2_{j,t})  \notag \\
&+\sum\limits_{j=1}^{N_p}(1-r^t_j)\hat{z}^t_j(m^3_{j,t}-m^4_{j,t})) \notag \\
&-\sum\limits_{i \in \mathcal{N}_g}(RU^1_i+p^0_i)\beta^1_{i,1} -\sum\limits_{i \in \mathcal{N}_g,t\in \mathcal{T}\setminus\{1\}}RU^t_i\beta^1_{i,t}\notag \\
&-\sum\limits_{i \in \mathcal{N}_g}(RD^1_i-p^0_i)\beta^2_{i,1} -\sum\limits_{i \in \mathcal{N}_g,t\in \mathcal{T}\setminus\{1\}}RD^t_i\beta^2_{i,t}\notag  \\
&+\sum\limits_{i \in \mathcal{N}_g, t\in \mathcal{T}}(-p^{\max}_i\kappa^1_{i,t}+p^{\min}_i\kappa^2_{i,t}-p^{\max}_ix^t_i\lambda^1_{i,t} \notag \\ &+p^{\min}_ix^t_i\lambda^2_{i,t}) \notag \\
 &+ \sum\limits_{t \in \mathcal{T}}\bar{q}^t \iota_{t} - \sum\limits_{i \in \mathcal{N}_g, t \in \mathcal{T}}\bar{q}^t_i\omega_{i,t} \\
\mbox{subject to } &\alpha^+_t-\alpha^-_t -\beta^1_{i,t}+\beta^1_{i,t+1} +\beta^2_{i,t}-\beta^2_{i,t+1} \notag \\
 &-\kappa^1_{i,t}+\kappa^2_{i,t}  -\lambda^1_{i,t}+\lambda^2_{i,t} \le C^t_i, \qquad \forall i \in \mathcal{N}_g, t\in \mathcal{T} \\
& -\kappa^1_{i,t}+\kappa^2_{i,t}+\iota_{t}-\omega_{i,t} \le 0, \qquad \forall i \in \mathcal{N}_g, t\in \mathcal{T} \\
& m^1_{j,t} \le \theta_m\zeta^t_j, \qquad \forall j \in \mathcal{N}_d,  t\in \mathcal{T} \\
& m^1_{j,t} \le \alpha_t^+, \qquad \forall j \in \mathcal{N}_d,  t\in \mathcal{T} \\
& m^1_{j,t} \ge \alpha_t^+ - \theta_m(1-\zeta^t_j), \qquad \forall j \in \mathcal{N}_d,  t\in \mathcal{T} \\
& m^2_{j,t} \le \theta_m\zeta^t_j, \qquad \forall j \in \mathcal{N}_d,  t\in \mathcal{T} \\
& m^2_{j,t} \le \alpha_t^-, \qquad \forall j \in \mathcal{N}_d,  t\in \mathcal{T} \\
& m^2_{j,t} \ge \alpha_t^- - \theta_m(1-\zeta^t_j), \qquad \forall j \in \mathcal{N}_d,  t\in \mathcal{T} \\
& m^3_{j,t} \le \theta_m\eta^t_j, \qquad \forall j \in \mathcal{N}_p,  t\in \mathcal{T} \\
& m^3_{j,t} \le \alpha_t^+, \qquad \forall j \in \mathcal{N}_p,  t\in \mathcal{T} \\
& m^3_{j,t} \ge \alpha_t^+ - \theta_m(1-\eta^t_j), \qquad \forall j \in \mathcal{N}_p,  t\in \mathcal{T} \\
& m^4_{j,t} \le \theta_m\eta^t_j, \qquad \forall j \in \mathcal{N}_p,  t\in \mathcal{T} \\
& m^4_{j,t} \le \alpha_t^-, \qquad \forall j \in \mathcal{N}_p,  t\in \mathcal{T} \\
& m^4_{j,t} \ge \alpha_t^- - \theta_m(1-\eta^t_j), \qquad \forall j \in \mathcal{N}_p,  t\in \mathcal{T} \\
\end{flalign}
\begin{flalign}
& \frac{1}{N_p} \left( \sum\limits_{k = 1}^{N_p}\sum\limits_{t = 1}^{\mathcal{T}}((\bar{z}^t_k- \eta^t_k \hat{z}^t_k)(1-r^t_k)) \right) \notag\\
 &- \sum_{t=1}^{\mathcal{T}} ((\bar{z}^t_j- \eta^t_j \hat{z}^t_j)(1-r^t_j)) \notag \\ &= a^+_j - a^-_j  \ \ \ (\forall j  \in \mathcal{N}_p) \\
& \beta^1,\beta^2,\kappa^1,\kappa^2,\lambda^1,\lambda^2,\iota,\omega,\alpha^+,\alpha^-, \notag\\ &m^1,m^2,m^3,m^4, a^+_j, a^-_j \ge 0 \\
& \zeta, \eta, \delta \  \in \{0,1\}
\end{flalign}
$\beta^1,\beta^2,\kappa^1,\kappa^2,\lambda^1,\lambda^2,\iota,\omega,\alpha^+,\alpha^-$ are dual variables.
$m^1,m^2,m^3,$ $m^4$ are new variables to handle bilinear term.
$a^+_j, a^-_j$ are also new variables for representing absolute value function constraints.

\section{Robust Problem}\label{Appendix:Robust Problem}
\begin{flalign}
\min\limits_{x,u,v,r,w} & \sum\limits_{t=1}^T\sum\limits_{i=1}^{N_g}x^t_i F^t_i+u^t_i S^t_i+v^t_i G^t_i+\sum\limits_{j=1}^{N_p} \Pi^t_j(r^t_j)+w\\
\mbox{subject to }& x^{t-1}_i-x^t_i +u^t_i  \ge 0,\qquad \forall i \in \mathcal{N}_g,t\in \mathcal{T} \\
& x^{t}_i-x^{t-1}_i +v^t_i  \ge 0,\qquad \forall i \in \mathcal{N}_g,t\in \mathcal{T} \\
& x^{t}_i-x^{t-1}_i \le x^{\tau}_i \qquad \forall \tau \in [t+1,\min\{t \notag \\ &\hspace{10mm} +\mbox{MinUp}_i-1,T\}],t\in[2,T]\\
& x^{t-1}_i-x^{t}_i \le 1-x^{\tau}_i \qquad \forall \tau \in [t+1,\min\{t \notag \\ &\hspace{10mm} +\mbox{MinDw}_i-1,T\}],t\in[2,T]\\
& w \ge -\chi L^1(a^+s, a^-s)+F(x,u,v,r)^\top y^s, \ \ \forall y^s \in \mathcal{S} 
\end{flalign}

\end{document}